\documentclass{article}
\date{}
\pdfoutput=1
\usepackage{hyperref}
\usepackage{indentfirst}
\usepackage{amsmath}
\usepackage{amssymb}
\title{Spherically symmetric Finsler metrics in $R^n$}
\author{Linfeng Zhou}
\begin{document}
\maketitle
\begin{center}
\textit{Department of Mathematics}\\
\textit{East China Normal University, Shanghai, 200241, China} \\
\textit{E-mail: lfzhou@math.ecnu.edu.cn}
\end{center}

\begin{abstract} In this paper, we give the general form of
spherically symmetric Finsler metrics in $R^n$ and surprisedly find
that many well-known Finsler metrics belong to this class. Then we
explicitly express projective metrics of this type. The necessary
and sufficient conditions that projective Finsler metrics with
spherical symmetry have constant flag curvature are also obtained.
\\

\end{abstract}

\newtheorem{Th}{Theorem}[section]
\newtheorem{Prop}[Th]{Propositon}
\newtheorem{Col}[Th]{Corollary}
\newtheorem{Lem}[Th]{Lemma}
\newtheorem{Ex}[Th]{Example}
\newtheorem{Con}[Th]{Conjecture}

\section{Introduction}

When studying Finsler geometry, one often encounters the intricacy
of calculation. So if we consider the Finsler metrics with a certain
symmetry, it may make things much easier.
 In general relativity, when
looking for a solution of the vacuum Einstein field equations
describing the gravitational field which is spherically symmetric,
we will obtain the Schwarzschild solution in four dimensional
time-space \cite{Yu}. In the process, the condition of spherical
symmetry plays a very important role. It can greatly simplify the
computation.
 Motivated by this idea, we investigate
spherically symmetric Finsler metrics in $R^n$ in this paper.

Similarly with the definition in general relativity, a spherically
symmetric Finsler metric means that it is invariant under any
rotations in $R^n$. In another words, the vector fields generated by
rotations are the Killing fields of the Finsler metric. Therefore we
firstly introduce the Killing field equation in Finsler geometry,
which generalizes the Killing field equation in Riemannian case
\cite{Pe}.

By solving the equations of Killing fields generated by rotations,
we firstly determine the structure of spherically symmetric Finsler
metrics $F(x,y)$ in $R^n$ in theorem \ref{th0}: $F$ must have the
form $F=\phi(|x|,|y|,\langle x,y\rangle)$. With a little surprise,
many well-known examples including Bryant metric \cite{Br} belong to
this type. Furthermore, spherically symmetric Finsler metrics are
not always $(\alpha,\beta)$ metrics. It is quite valuable to study
this type of Finsler metrics.

To character projective Finsler metrics in $R^n$ is a very important
problem. It is connected with Hilbert's Fourth problem \cite{Sh1}.
Thus next we discuss those spherically symmetric Finsler metrics
$F(x,y)$ in $R^n$ which are projective. We express them explicitly
by $F=\int f(\frac{v^2}{u^2}-r^2)du+g(r)v$ in theorem \ref{th1} by
using famous Rapcs\'{a}k's lemma \cite{Ra}. As we know, a projective
Finsler metric is of scalar curvature. It is natural to ask which
metrics have constant flag curvature among the projective
spherically symmetric Finsler metrics. We obtain a sufficient and
necessary condition in theorem \ref{th2}, which is two partial
differential equations. From these equations, perhaps one can find
some new examples of Finsler metrics with constant flag curvature.

 Finally, we suggest a conjecture \ref{con1} by observing the examples as far as we know.

\section{The Killing field equation in Finsler geometry}
Suppose $F$ is a Finsler metric on an n-dimensional $C^{\infty}$
manifold $M$. Let $\Phi:M\rightarrow M$ be a diffeomorphism and
$\Phi_{\ast}: T_{x}M\rightarrow T_{\Phi(x)}M$ be the tangent map at
point $x$. $\Phi$ is called isometric map if it satisfies
\[F(\Phi(x),\Phi_{\ast}(y))=F(x,y)\]
where $y\in T_{x}M$.

A vector field $X$ on $M$ is called a Killing field if the
1-parameter group $\Phi_{t}$ generated by $X$ are isometric. As we
know in Riemann geometry, there is the following Killing field
equation \cite{Pe}
\[\mathcal{L}_{X}g=0.\]
Here $\mathcal{L}_{X}g$ denotes the Lie derivative of Riemannian
metric tensor $g$ on $M$. In local coordinate
$\{x^i,\frac{\partial}{\partial x^i}\}$, above equation can be
written as
\[\frac{\partial{g_{ij}}}{\partial x^p}X^p+g_{pj}\frac{\partial{X^p}}{\partial{x^i}}+g_{ip}\frac{\partial{X^p}}{\partial{x^j}}=0\]
where $g=g_{ij}dx^i\otimes dx^j$ and
$X=X^i\frac{\partial}{\partial{x^i}}$. Hence there must be a similar
Killing field equation in Finsler geometry. In fact, we have the
following theorem:
\begin{Th} \label {th00} Let $(M^n,F)$ be an n-dimensional smooth Finsler manifold. A
vector field $X$ is a Killing field on $M$. Then $X$ satisfies the
equation
\[\frac{\partial{g_{ij}}}{\partial x^p}X^p+g_{pj}\frac{\partial{X^p}}{\partial{x^i}}+g_{ip}\frac{\partial{X^p}}{\partial{x^j}}+2C_{ijp}\frac{\partial{X^p}}{\partial{x^k}}y^k=0\]
where
$g_{ij}=\frac{1}{2}\frac{\partial^2{F^2}}{\partial{y^i}\partial{y^j}}$
is coefficient of fundamental tensor,
$C_{ijp}=\frac{1}{2}\frac{\partial g_{ij}}{\partial{y^p}}$ is
coefficient of Cartan tensor and $X=X^i\frac{\partial}{\partial
{x^i}}$ under the local coordinate.
\end{Th}
\emph{Proof.} Suppose $\Phi_{t}$ is the 1-parameter group of $X$.
According to the definition of the Killing field in Finsler
geometry, we have
\[F(\Phi_{t}(x),(\Phi_{t})_{\ast}(y))=F(x,y).\]
Under the local coordinate that means
\[g_{ij}(\Phi_{t}(x),(\Phi_{t})_{\ast}(y))\frac{\partial{\Phi_t^i}}{\partial{x^k}}y^k\frac{\partial{\Phi_t^j}}{\partial{x^l}}y^l=g_{ij}(x,y)y^iy^j.\]
Notice that $X_x=\frac{d\Phi_{x}(t)}{dt}|_{t=0}$ and $\Phi_{0}$ is
identity. So taking derivative with respect to $t$ in above equation
and set $t=0$, we get
\[\frac{\partial{g_{ij}}}{\partial{x^p}}X^py^iy^j+\frac{\partial{g_{ij}}}{\partial{y^p}}\frac{\partial{X^p}}{\partial x^k}y^ky^iy^j
+g_{ij}y^p\frac{\partial{X^i}}{\partial{x^p}}y^j+g_{ij}y^i\frac{\partial{X^j}}{\partial{x^k}}y^k=0.\]
It is equivalent to
\[\frac{\partial{g_{ij}}}{\partial{x^p}}X^p+\frac{\partial{g_{ij}}}{\partial{y^p}}\frac{\partial{X^p}}{\partial x^k}y^k
+g_{pj}\frac{\partial{X^p}}{\partial{x^i}}+g_{ip}\frac{\partial{X^p}}{\partial{x^j}}=0.\]
By the definition of Cartan torsion, we obtain the result
immediately.\quad Q.E.D.
\\

\noindent\emph{Remark.} In the case of Riemannian metric, the Cartan
torsion vanishes. Hence the equation in theorem \ref{th00} coincides
with $\mathcal{L}_{X}g=0$.

\begin{Col}\label{col1} Let $(M^n,F)$ be an n-dimensional smooth Finsler manifold. A
vector field $X$ is a Killing field on $M$. Then $X$ satisfies
\[\frac{\partial{F}}{\partial{x^i}}X^i+\frac{\partial{F}}{\partial{y^i}}\frac{\partial{X^i}}{\partial{x^j}}y^j=0.\]
\end{Col}
\emph{Proof.} The conclusion can be obtained by contracting the
equation in theorem \ref{th00} with $y^i$ and $y^j$.\quad Q.E.D.

\section{The general form of spherically symmetric Finsler metrics in $R^n$}
We denote $\Omega$ a convex domain in $R^n$ and $F$ a Finsler metric
on $\Omega$. $(\Omega, F)$ is called spherically symmetric if
orthogonal matrix $O(n)$ is isometric map of $(\Omega, F)$. This
definition is equivalent to say that $(\Omega, F)$ is invariant
under any rotations in $R^n$. Hence there is a natural question:
what is the restriction on the metric $F$ if $F$ has spherical
symmetry? We give the following answer:
\begin{Th} \label{th0}
Let $F(x,y)$ be a Finsler metric on a convex domain $\Omega
\subseteq R^n$. $F(x,y)$ is spherically symmetric if and only if
there exists a positive function $\phi(r,u,v)$ s.t.
\[F(x,y)=\phi(|x|,|y|,\langle x,y\rangle)\]
where $|x|=\sqrt{(x^1)^2+\cdots +(x^n)^2}$,
$|y|=\sqrt{(y^1)^2+\cdots +(y^n)^2}$ and $\langle
x,y\rangle=x^1y^1+\cdots+x^ny^n$.
\end{Th}
\emph{Proof.} Suppose $F(x,y)$ is spherically symmetric on $\Omega
\subseteq R^n$. Choose $\{e_1,\cdots,e_n\}$ is the standard
orthonormal base in $R^n$ and denote $X_iOX_j$ the coordinate plane
spanned by $\{e_{i},e_{j}\}$. Consider a family rotations $\theta_t$
on coordinate plane $X_iOX_j$:
\[\begin{split}&\theta_t(x^1,\cdots,x^i,\cdots,x^j,\cdots,x^n)\\
=&(x^1,\cdots,x^i\cos t+x^j\sin t,\cdots,-x^i\sin t+x^j\cos
t,\cdots,x^n).\end{split}\] Obviously $\theta_t$ is a 1-parameter
group and isometric. So a Killing vector field X generated by
$\theta_t$ is
\[X=x^j\frac{\partial}{\partial x^i}-x^i\frac{\partial}{\partial x^j}.\]
By corollary \ref{col1}, we have the following equation
\begin{equation}\label{eq1}\frac{\partial F}{\partial x^i}x^j-\frac{\partial F}{\partial
x^j}x^i+\frac{\partial F}{\partial y^i}y^j-\frac{\partial
F}{\partial y^j}y^i=0.\end{equation} This equation is a first order
linear partial differential equation. It's characteristic equation
is given by
\[\frac{dx^i}{x^j}=-\frac{dx^j}{x^i}=\frac{dy^i}{y^j}=-\frac{dy^j}{y^i}.\]
Thus
\[(x^i)^2+(x^j)^2=c_1,\quad (y^i)^2+(y^j)^2=c_2,\quad x^iy^i+x^jy^j=c_3\]
are three independent first integrals. Hence the solution of
equation (\ref{eq1}) is
\[F=\overline{\phi}(x^1,\cdots,\widehat{x^i},\cdots,\widehat{x^j},\cdots,x^n,(x^i)^2+(x^j)^2,(y^i)^2+(y^j)^2,x^iy^i+x^jy^j).\]
Here $\widehat{x^i}$ means omitting the variable $x^i$.
 For $i,j$ are arbitrary numbers
from $1$ to $n$, so there are $\frac{n(n-1)}{2}$ Killing field
equations like (\ref{eq1}). Therefore, $F$ must have the following
form
\[\begin{split}F(x,y)&=\widetilde{{\phi}}((x^1)^2+\cdots+(x^n)^2,(y^1)^2+\cdots+(x^n)^2,x^1y^1+\cdots+x^ny^n)\\
&=\phi(|x|,|y|,\langle x,y\rangle).\end{split}\] The converse is
obvious by a direct computation.\quad Q.E.D.\\

Let $F=\phi(|x|,|y|,\langle x,y\rangle)$ where $\phi(r,u,v)$ is a
positive $C^{\infty}$ function with homogeneous of degree one with
respect to variable $u$ and $v$, let us find the condition for the
positivity of $(g_{ij}):=(\frac{1}{2}\frac{\partial ^2F^2}{\partial
y^i\partial y^j})$. It is easy to compute $g_{ij}$:
\[g_{ij}=\frac{\phi\phi_u}{u}\delta_{ij}+(\phi_v^2+\phi\phi_{vv})x^ix^j+(\frac{\phi_u^2+\phi\phi_{uu}}{u^2}-\frac{\phi\phi_u}{u^3})y^iy^j+(\frac{\phi_u\phi_v+\phi\phi_{uv}}{u})(x^iy^j+x^jy^i)\]
where $u:=|y|$. Thus we can obtain \cite{CS}
\[det(g_{ij})=(\frac{\phi}{u})^{n+1}\phi_u^{n-2}[\phi_u+(|x|^2|y|^2-\langle
x,y\rangle^2)\frac{\phi_{vv}}{u}].\]

\begin{Lem}
Suppose a positive $C^{\infty}$ function $\phi(r,u,v)$ is
homogeneous of degree one with respect to $u$ and $v$. If $\phi$
satisfies that
\[\phi_u>0,\phi_{uu}\geq 0 \]
when $u>0$ and $r\ge 0$, then $F=\phi(|x|,|y|,\langle x,y\rangle)$
is a Finsler metric.
\end{Lem}
\emph{Proof.} Since $\phi(r,u,v)$ is homogeneous of degree one, we
have
\[\phi_{uu}u+\phi_{uv}v=0.\]
So
\[\phi_{uu}=-\frac{v}{u}\phi_{uv}.\]
Similarly we have
\[\phi_{uv}=-\frac{v}{u}\phi_{vv}.\]
Thus we obtain that
\[\phi_{uu}=(\frac{v}{u})^2\phi_{vv}.\]
Hence the condition $\phi_u>0,\phi_{uu}\geq 0$ is equivalent to
\[\phi_u>0,\phi_{vv}\geq 0.\]
By the formula of $det(g_{ij})$ computed above, one can easily see
that the matrix $(g_{ij})$ is positive.\quad Q.E.D.\\

In fact, many classical Finsler metrics are spherically symmetric
\cite{BCS} \cite{CS}.
\begin{Ex}[Klein model] Let $B^n \subset R^n$ be the standard unit
ball and let
\[\alpha(x,y):=\frac{\sqrt{|y|^2-(|x|^2|y|^2-\langle x,y\rangle^2)}}{1-|x|^2},\quad y\in T_xB^n.\]
$\alpha(x,y)$ is a Riemannnian metric on $B^n$. It is projective and
has constant flag curvature $K=-1$.
\end{Ex}

\begin{Ex}[Funk metric]
A Randers metric $F$ is defined on the standard unite ball $B^n$:
\[F(x,y):=\frac{\sqrt{|y|^2-(|x|^2|y|^2-\langle x,y\rangle^2)}+\langle
x,y\rangle}{1-|x|^2}.\] It is also projective and has constant flag
curvature $K=-\frac{1}{4}$.
\end{Ex}

\begin{Ex}[Berwald metric]
An $(\alpha,\beta)$ metric $F$ is also defined on $B^n$:
\[F(x,y):=\frac{(\sqrt{|y|^2-(|x|^2|y|^2-\langle x,y\rangle^2)}+\langle x,y\rangle)^2}{(1-|x|^2)^2\sqrt{|y|^2-(|x|^2|y|^2-\langle x,y\rangle^2)}}.\]
$F$ is projective and has constant flag curvature $K=0$.
\end{Ex}

\begin{Ex} [projective spherical model]
Let $S^n \subset R^{n+1}$ be the standard unit sphere. The standard
inner product $\langle,\rangle$ in $R^{n+1}$ induced a Riemannian
metric on $S^n$: for $x\in S^n$, let
\[\alpha:=|y|,\quad y\in T_xS^n\subset R^{n+1}.\]
Let $S^n_{+}$ denote the upper hemisphere and let
$\psi_{+}:R^n\rightarrow S_{+}^n$ be the projection map defined by
\[\psi_{+}(x):=(\frac{x}{\sqrt{1+|x|^2}},\frac{1}{1+|x|^2}).\]
The pull-back metric on $R^n$ from $S^n_{+}$ by $\psi_{+}$ is given
by
\[\alpha(x,y):=\frac{\sqrt{|y|^2+(|x|^2|y|^2-\langle x,y \rangle^2)}}{{1+|x|^2}},\quad y\in T_xR^n\]
$(R^n,\alpha(x,y))$ is projectively flat and has constant flag
curvature $K=1$.
\end{Ex}

\begin{Ex}[Bryant metric]
Denote
\[\begin{split}
A:=&(\cos(2\alpha)|y|^2+(|x|^2|y|^2-\langle x,y\rangle^2))^2+(\sin(2\alpha)|y|^2)^2,\\
B:=&\cos(2\alpha)|y|^2+(|x|^2|y|^2-\langle x,y\rangle^2),\\
C:=&\sin(2\alpha)\langle x,y\rangle,\\
D:=&|x|^4+2\cos(2\alpha)|x|^2+1.
\end{split}\]
For an angle $\alpha$ with $0\leq\alpha<\frac{\pi}{2}$, Bryant
metric $F$ is defined by
\[F:=\sqrt{\frac{\sqrt{A}+B}{2D}+(\frac{C}{D})^2}+\frac{C}{D}\]
on the whole region $R^n$. As we know it is projective and has
constant flag curvature $K=1$.
\end{Ex}

From above examples, we can see that spherically symmetric Finsler
metrics don't always belong to $(\alpha,\beta)$ metrics. So it is
meaningful to study projective metrics of this type with constant
flag curvature \cite{LS}.

\section{Projective spherically symmetric Finsler metrics in $R^n$}

A Finlser metric $F$ in $R^n$ is called projective metric, if its
geodesics are straight lines. Since spherically symmetric Finsler
metrics have very nice symmetry, when discussing projective metrics
of this type in $R^n$, we can obtain a quite simple result compared
with projectively flat $(\alpha,\beta)$ metrics \cite{Sh2}. Before
stating our result, we need an important lemma about projective
Finsler metrics.

\begin{Lem} (Rapcs\'{a}k \cite{Ra}) \label{lem1}Let $F(x,y)$ be a Finsler metric on an
open subset $\mathcal{U}\in R^n$. $F(x,y)$ is projective on
$\mathcal{U}$ if and only if it satisfies
\[F_{x^ky^l}y^k=F_{x^l}.\]
In this case, the projective factor $P(x,y)$ is given by
\[P=\frac{F_{x^k}y^k}{2F}.\]
\end{Lem}

\begin{Th} \label{th1}Suppose $F$ is a spherically symmetric Finsler
metric on a convex domain $\Omega\in R^n$, $F$ is projective if and
only if there exist smooth functions $f(t)>0$ and $g(r)$ s.t.
\[\phi(r,u,v)=\int f(\frac{v^2}{u^2}-r^2)du+g(r)v\]
where $F(x,y)=\phi(|x|,|y|,\langle x,y\rangle)$.
\end{Th}
\emph{Proof.} By lemma \ref{lem1}, $F$ is projective if and only if
$F$ satisfies
\begin{equation}\label{eq2}F_{x^l}=F_{y^lx^k}y^k.\end{equation}
If $F$ is spherically symmetric, then there exists $\phi(r,u,v)$
s.t.
\[F(x,y)=\phi(|x|,|y|,\langle x,y\rangle).\]
So
\[F_{x^l}=\phi_r\frac{x^l}{|x|}+\phi_vy^l=\phi_r\frac{x^l}{r}+\phi_vy^l\]
and
\[\begin{split}F_{y^lx^k}y^k&=(\phi_{y^l})_{x^k}y^k=(\phi_u\frac{y^l}{|y|}+\phi_{v}x^l)_{x^k}y^k\\
&=\sum_{k=1}^n(\phi_{ru}\frac{x^k}{|x|}\frac{y^l}{|y|}+\phi_{uv}y^k\frac{y^l}{|y|}+\phi_{rv}\frac{x^k}{|x|}x^l+\phi_{vv}y^kx^l+\phi_v\delta^{l}_k)y^k\\
&=\phi_{ru}\frac{\langle
x,y\rangle}{|x||y|}y^l+\phi_{uv}|y|y^l+\phi_{rv}\frac{\langle
x,y\rangle}{|x|}x^l+\phi_{vv}|y|^2x^l+\phi_vy^l\\
&=(\phi_{rv}\frac{v}{r}+\phi_{vv}u^2)x^l+(\phi_{ru}\frac{v}{ru}+\phi_{uv}u+\phi_v)y^l\end{split}\]
where $r:=|x|$, $u:=|y|$, $v:=\langle x,y\rangle$.
 Thus (\ref{eq2}) holds if and only if $\phi$ satisfies
\begin{equation}\label{eq3}
\left\{ \begin{aligned}
         \phi_{rv}\frac{v}{r}+\phi_{vv}u^2=\frac{\phi_r}{r} \\
                 \phi_{uv}u+\phi_{ru}\frac{v}{ru}=0 .
          \end{aligned} \right.\end{equation}
In fact, the two equations in (\ref{eq3}) are equivalent by only
noticing that
\[\phi_{rv}=\frac{\phi_r-u\phi_{ru}}{v},\quad \phi_{vv}=-\frac{u\phi_{uv}}{v}\]
for $\phi$ being homogeneous of degree one. Hence $F$ is projective
if and only if $\phi$ satisfies
\begin{equation}\label{eq4}\phi_{uv}u+\phi_{ru}\frac{v}{ru}=0.\end{equation} This is a first order
linear partial differential equation with respect to $\phi_u$. It's
characteristic equation is
\[\frac{dv}{u}=\frac{dr}{\frac{v}{ru}}.\]
So
\[\phi_u=\tilde{f}(v^2-u^2r^2,u^2)=\tilde{f}(\frac{v^2}{u^2}-r^2,1)=f(\frac{v^2}{u^2}-r^2)\]
is the solution of equation (\ref{eq4}). Here it should be pointed
out that in above equation, we use the homogeneous property of
$\phi_u$. Thus there exists a function $c(r,v)$ s.t.
\[\phi=\int f(\frac{v^2}{u^2}-r^2)du+c(r,v).\]
Again noticing the homogeneity of $\phi$, we conclude that
\[c(r,v)=g(r)v.\]
Thus complete the proof. \quad Q.E.D.\\

As we know, projective Finsler metric is of scalar curvature. So
from above theorem, we can find many spherically symmetric Finsler
metrics having scalar curvature. Now let us study those projective
Finsler metrics which have constant flag curvature among this type
in $R^n$. We need a lemma first.
\begin{Lem} \cite{Sh1} \label{lem2}
Suppose $F=F(x,y)$ is a projective Finlser metric on a convex domain
$\Omega \subseteq R^n$. Then $F$ has constant flag curvature
$K=\lambda$ if and only if projective factor $P$ satisfies
\[P_{x^k}=PP_{y^k}-\lambda FF_{y^k}\]
where $P:=\frac{F_{x^m}y^m}{2F}$.
\end{Lem}

With this lemma, we have the following conclusion:
\begin{Th}\label{th2} Let a spherically symmetric Finsler metric $F=\phi(|x|,|y|,\langle
x,y\rangle)$ be projective on a convex domain $\Omega \subseteq
R^n$. Then $F$ has constant flag curvature $K=\lambda$ if and only
if $\phi(r,u,v)$ satisfies

\begin{equation}\label{eq5}
\left\{ \begin{aligned}
          4\lambda r\phi^4\phi_u+r\phi_uQ^2-4ru\phi\phi_vQ+4u\phi^2\phi_r=0\\
          4\lambda r\phi^4\phi_v+r\phi_vQ^2+2\phi^2 Q_r-4\phi\phi_rQ=0
          \end{aligned} \right.\end{equation}
where $Q:=\frac{v}{r}\phi_r+u^2\phi_v$.
\end{Th}
\emph{Proof.} From Lemma \ref{lem2}, $F$ has constant flag curvature
$K=\lambda$ if and only if
\begin{equation}\label{eq6}
P_{x^k}=PP_{y^k}-\lambda FF_{y^k}
\end{equation}
where projective factor $P=\frac{F_{x^m}y^m}{2F}$. Now
$F=\phi(r,u,v)$, $r=|x|$, $u=|y|$ and $v=\langle x,y\rangle$. Hence
\[\begin{split}
P=&\frac{1}{2\phi}(\frac{v}{r}\phi_r+u^2\phi_v),\\
P_{x^k}=&-\frac{1}{2\phi^2}(\frac{\phi_r}{r}x^k+\phi_vy^k)(\frac{v}{r}\phi_r+u^2\phi_v)\\
&+\frac{1}{2\phi}(\phi_{rr}\frac{v}{r^2}x^k+\phi_{rv}\frac{v}{r}y^k+\frac{\phi_r}{r}y^k-\phi_r\frac{v}{r^3}x^k+\phi_{rv}\frac{u^2}{r}x^k+\phi_{vv}u^2y^k),\\
PP_{y^k}=&P(-\frac{1}{2\phi^2})(\frac{\phi_u}{u}y^k+\phi_vx^k)(\frac{v}{r}\phi_r+u^2\phi_v)\\
&+\frac{P}{2\phi}(\phi_{ru}\frac{v}{ru}y^k+\phi_{rv}\frac{v}{r}x^k+\frac{\phi_r}{r}x^k+u\phi_{uv}y^k+\phi_{vv}u^2x^k+2\phi_vy^k),\\
\lambda FF_{y^k}=&\lambda \phi(\frac{\phi_u}{u}y^k+\phi_{v}x^k).
\end{split}\]
Substituting above equations into (\ref{eq6}), we can know that
(\ref{eq6}) holds if and only if
\begin{equation}\label{eq7}\begin{split}
&-\frac{\phi_v}{2\phi^2}(\frac{v}{r}\phi_r+u^2\phi_v)+\frac{1}{2\phi}(\phi_{rv}\frac{v}{r}+\frac{\phi_r}{r}+\phi_{vv}u^2)\\
&=-\frac{P}{2\phi^2}\frac{\phi_u}{u}(\phi_r\frac{v}{r}+\phi_vu^2)+\frac{P}{2\phi}(\phi_{ru}\frac{v}{ru}+u\phi_{uv}+2\phi_v)-\lambda
\phi\frac{\phi_{u}}{u}
\end{split}\end{equation}
and
\begin{equation}\label{eq8}
\begin{split}
       &-\frac{\phi_r}{2\phi^2r}(\frac{v}{r}\phi_r+u^2\phi_v)+\frac{1}{2\phi}(\phi_{rr}\frac{v}{r^2}-\phi_r\frac{v}{r^3}+\phi_{rv}\frac{u^2}{r})\\
       &= -\frac{P}{2\phi^2}\phi_{v}(\frac{v}{r}\phi_r+u^2\phi_v)
         +\frac{P}{2\phi}(\phi_{rv}\frac{v}{r}+\frac{\phi_r}{r}+\phi_{vv}u^2)-\lambda \phi\phi_v. \end{split}
\end{equation}
Noticing that $F$ is projective, so $\phi$ satisfies equation
(\ref{eq3}) in theorem \ref{th1}
\[ \left\{ \begin{aligned}
         \phi_{rv}\frac{v}{r}+\phi_{vv}u^2=\frac{\phi_r}{r} \\
                 \phi_{uv}u+\phi_{ru}\frac{v}{ru}=0 .
          \end{aligned} \right.\]
Using above two equations and substituting the formula of $P$, we
can simplify equation (\ref{eq7}), (\ref{eq8}) and obtain the
result.\qquad Q.E.D.\\

\noindent\emph{Remark.} \begin{enumerate}
\item [1)]Although we have
known the general form of projective spherically symmetric metric,
it is difficult to solve the equations (\ref{eq5}) in theorem
\ref{th2} directly.
\item [2)] It can be verified that the examples in section 3 satisfy
our theorem \ref{th1} and theorem \ref{th2} by Maple.
\end{enumerate}

Finally we propose the following conjecture by observing the
examples in section 3:
\begin{Con} \label{con1}
Suppose $F=F(x,y)$ is a projective spherically symmetric Finlser
metric with constant flag curvature on a convex domain $\Omega
\subseteq R^n$. If $F$ is reversible, then $F$ must be Riemannian.
\end{Con}

 \LaTeX
\end{document}